\documentclass[letterpaper,final,twoside,reqno]{amsart}
\usepackage{xspace}
\usepackage{amssymb}

\newcommand{\baseRing}[1]{\ensuremath{\mathbb{#1}}}
\newcommand{\Z}{\baseRing{Z}}
\newcommand{\R}{\baseRing{R}}
\newcommand{\C}{\baseRing{C}}
\newcommand{\idM}{\ensuremath{Id}}
\newcommand{\pd}[2]{{\frac{\partial {#1}}{\partial {#2}}}}
\newcommand{\pdii}[3]{\frac{\partial^2 {#1}}{\partial z_{#2} \partial
    z_{#3}}}
\newcommand{\pdiii}[4]{\frac{\partial^3 {#1}}{\partial z_{#2} \partial
    z_{#3} \partial z_{#4}}}
\newcommand{\CB}{\ensuremath{\mathcal{B}}\xspace}
\newcommand{\jgg}{\ensuremath{\mathfrak{g}}\xspace}
\newcommand{\jgp}{\ensuremath{\mathfrak{p}}\xspace}
\newcommand{\liesl}{\ensuremath{\mathfrak{s}\mathfrak{l}}\xspace}
\newcommand{\jdef}[1]{\emph{#1}}
\newcommand{\h}{\hbar}

\theoremstyle{plain}
\newtheorem{theorem}{Theorem}[section]

\newtheorem{proposition}[theorem]{Proposition}
\newtheorem{lemma}[theorem]{Lemma}

\theoremstyle{definition}
\newtheorem{definition}[theorem]{Definition}
\newtheorem{remark}[theorem]{Remark}
\newtheorem{example}[theorem]{Example}

\numberwithin{equation}{section}

\DeclareMathOperator{\aut}{Aut}
\DeclareMathOperator{\gr}{Gr}
\DeclareMathOperator{\vspan}{Span}
\DeclareMathOperator{\res}{Res}

\DeclareMathOperator{\sym}{Sym}
\DeclareMathOperator{\jend}{End}

\newcommand{\Script}[1]{\ensuremath{{\mathcal{#1}}}}

\newcommand{\UU}{\Script{U}}
\newcommand{\FF}{\Script{F}}
\newcommand{\KK}{\Script{K}}
\newcommand{\QQ}{\Script{Q}}    
\newcommand{\VV}{\Script{V}}

\newcommand{\CC}{\Script{C}}

\newcommand{\ie}{\textsl{i.e.}\xspace}

\newcommand{\ti}[1]{\tilde{#1}}
\newcommand{\conj}{\overline}

\newcommand{\GC}{\ensuremath{G_\C}\xspace}
\newcommand{\GR}{\ensuremath{G_\R}\xspace}
\newcommand{\gl}{\ensuremath{\mathfrak{gl}}}
\newcommand{\lgr}{\ensuremath{\mathfrak{g}_\R}}
\newcommand{\DD}{\ensuremath{D}\xspace}
\newcommand{\DC}{\ensuremath{\check{D}}\xspace}

\begin{document}

\title{Frobenius modules and Hodge asymptotics}

\author{Eduardo Cattani} 
\address{Department of Mathematics and Statistics\\ University of
  Massachusetts\\ Amherst\\ MA 01003-9305\\USA.}
\email{cattani@math.umass.edu}
\thanks{E. Cattani was partially supported by
    NSF Grant DMS-0099707}

\author{Javier Fernandez}
\address{Department of Mathematics\\ University of Utah\\ Salt Lake
  City\\ UT 84112-0090\\USA.}
\email{jfernand@math.utah.edu}

\bibliographystyle{amsplain}



\begin{abstract}
  We exhibit a direct correspondence between the potential defining
  the $H^{1,1}$ small quantum module structure on the cohomology of a
  Calabi-Yau manifold and the asymptotic data of the $A$-model
  variation of Hodge structure.  This is done in the abstract context
  of polarized variations of Hodge structure and Frobenius modules.
\end{abstract}

\maketitle



\section{Introduction}
\label{sec:intro}

The even cohomology of a compact smooth manifold is a Frobenius algebra
with respect to the cup product and the intersection form.  For a
compact, K\"ahler manifold $X$, multiplication by a K\"ahler class
defines a representation of the Lie algebra $\liesl(2)$ on the full
cohomology $H^*(X,\C)$, whose semisimple element induces the standard
\Z-grading. This is the content of the Hard Lefschetz Theorem.
Beginning with the formulation of the Mirror Symmetry phenomenon
\cite{ar:CDGP-pair}, there has been considerable interest in studying
the simultaneous action on cohomology of the K\"ahler cone $\KK$ of
$X$.  Looijenga and Lunts~\cite{ar:LL-lefschetz_modules} have shown
that the copies of $\liesl(2)$ associated with the elements of $\KK$
generate a semisimple Lie algebra and have studied some of their
properties.  Another point of view, introduced in~\cite{ar:CKS2},
consists in studying $H^*(X,\C)$ as a mixed Hodge structure which
splits over $\R$ and is polarized by the action of every K\"ahler
class.  Hence, the crucial information is contained in the structure
of $H^*(X,\C)$ as a $\sym H^{1,1}$-module.  In particular, it follows
from~\cite[Prop.  4.66]{ar:CKS} that we may define a polarized
variation of Hodge structure on $H^*(X,\C)$ parametrized by the
complexified K\"ahler cone of $X$.  If a polyhedral cone of K\"ahler
classes is chosen, this variation becomes a nilpotent orbit in the
sense of Schmid~\cite{ar:S-vhs}.  This approach has proved fruitful in
the study of mixed Lefschetz theorems~\cite{ar:CKS2}.

Quantum cohomology is a deformation of the cup product on $H^*(X,\C)$
defined in terms of the Gromov-Witten potential ---a generating
function for certain enumerative invariants.  If $X$ is a Calabi-Yau
manifold, the action of $H^{1,1}$ on $\oplus H^{p,p}(X)$, with respect
to the small quantum product, leads to a variation of Hodge structure,
called the $A$-model variation by Morrison~\cite{ar:morrison-aspects}.
A local variation of Hodge structure is described by an algebraic
component ---the nilpotent orbit--- and an analytic part described by
a holomorphic map with values in a graded component of a nilpotent Lie
algebra.  For the $A$-model variation the nilpotent orbit is the one
described in the previous paragraph.

Both Frobenius algebras and polarized variations of Hodge structure
have been extensively studied in the recent physics literature.
Variations of Hodge structure appear, for instance, in connection with
the tree level amplitudes of twisted $N=2$ theories --the $B$-model--
and, for Calabi-Yau threefolds, as special geometry
(\cite{ar:BCOV-Kodaira_Spencer,ar:Cecotti-special,ar:CV-top_anti-top}).
On the other hand, $2D$ topological field theories are equivalent to
Frobenius algebras. Families of these algebras were also considered:
the tangent bundle of the moduli space of topological conformal field
theories has, on each fiber, a Frobenius algebra structure
(\cite{th:Dijkgraaf,ar:DVV-topological_strings}). A relation between
the two objects arises in mirror symmetry via the equivalence of the
$A$ and $B$ model correlation functions
(\cite{ar:CDGP-pair,ar:Greene-strings_CY,bo:CK-mirror,ar:morrison-aspects}).
What is perhaps not so well known is a direct construction due to
Morrison of a variation of Hodge structure based on the
$A$-model~\cite{ar:morrison-aspects}. In this paper we show a
correspondence between any polarized variation of Hodge structure with
appropriate degenerating behavior and a certain sub-structure of a
family of Frobenius algebras.  Our main result is to exhibit a simple,
direct correspondence between the holomorphic data of the variation
and the (small) quantum potential in such a way that the horizontality
equation of a variation of Hodge structure corresponds to a graded
component of the WDVV equations.

We will work throughout in the setting of abstract variations of Hodge
structure.  The analogous abstract notion on the ``quantum'' side is
that of a Frobenius module introduced in
Section~\ref{sec:frobenius_modules} and their deformations defined by
potentials encoding the essential properties of a graded portion of
the Gromov-Witten potential.

The paper is organized as follows.  In \S\ref{sec:hodge_theory} we
review the asymptotic description of variations.
Theorem~\ref{th:improved_2.8} contains the algebraic and analytic
characterization of local variations.  We also recall the notion of
maximally unipotent boundary points and of canonical coordinates
\cite{ar:CDGP-pair,ar:morr-picard-fuchs,ar:Del-local_behavior}.  In
Section~\ref{sec:frobenius_modules} we define Frobenius modules and
their deformations.  Section~\ref{sec:correspondence} is devoted to
the proof of our main result, Theorem~\ref{th:correspondence_grl},
which establishes an equivalence between local variations with
appropriate behavior at the boundary and quantum potentials.  Finally,
in \S\ref{sec:a_model_variation} we review the construction of the
$A$-model variation and show that it coincides with the one
constructed in Theorem~\ref{th:correspondence_grl}.  As a byproduct,
we obtain a direct proof that the $A$-model variation is indeed a
polarized variation of Hodge structure.

We note that the $A$-model variation involves only the small quantum
module structure.  In the case of Hodge structures of weights $3$, $4$
and $5$, corresponding to threefolds, fourfolds and fivefolds, the
module structure suffices to recover the full quantum algebra, so that
our results extend the previously known correspondences
(\cite{ar:greg-higgs,ar:CF-asymptotics}) in weights $3$ and $4$.
Also, the full quantum algebra can be recovered if it is assumed to be
generated, in the geometric context, by $H^{1,1}$. In this last case,
the family of Frobenius algebras obtained from a variation of Hodge
structure can be seen as a Frobenius manifold. These matters will be
analyzed elsewhere~\cite{ar:FP-opposite}.  S.
Barannikov~\cite{ar:Ba-projective,ar:Ba-quantum_periods,ar:Ba-generalized_periods}
has introduced the notion of semi-infinite variations of Hodge
structure to deal with the full quantum algebra.  He has also shown
that, for projective complete intersections, the $A$-model variation
is of geometric origin and coincides with the polarized variation of
Hodge structure of the mirror family.

Finally, we wish to thank Gregory Pearlstein for his very helpful comments.



\section{Hodge theory preliminaries}
\label{sec:hodge_theory}

In this section we briefly review the asymptotic description of
variations of Hodge structure.  We refer to
\cite{ar:Gri-periods-1,ar:S-vhs,ar:CK-luminy,ar:CF-asymptotics} for
details and proofs.

A (real) \jdef{variation of Hodge structure} $\VV$ over a connected
complex manifold $M$ consists of a holomorphic vector bundle $\VV\to
M$, a flat connection $\nabla$ on $\VV$ with quasi-unipotent
monodromy, a flat real form $\VV_{\R} \subset \VV$, and a finite
decreasing filtration $\FF$ of $\VV$ by holomorphic subbundles ---the
\jdef{Hodge filtration}--- satisfying
\begin{eqnarray}\label{eq:horizontality}
  \nabla\FF^p & \subset & \Omega^1_M
  \otimes \FF^{p-1}\quad\hbox{(Griffiths' Transversality) and}
  \\ \label{eq:opposed}
  \VV& = & \FF^p \oplus \conj{\FF}^{k-p+1}
\end{eqnarray}
for some integer $k$ ---the \jdef{weight} of the variation--- and
where barring denotes conjugation relative to $\VV_{\R}$.  As a
$C^{\infty}$-bundle, $\VV$ may then be written as a direct sum
\begin{equation}\label{bundledecomposition}
  \VV = \bigoplus_{p+q=k}\ \VV^{p,q}\ ,\quad \quad \VV^{p,q} := \FF^p
  \cap \conj{\FF}^q\,;
\end{equation}
the integers $h^{p,q} := \dim\,\VV^{p,q}$ are the \jdef{Hodge numbers}.
A \jdef{polarization} of the variation is a flat non-degenerate
bilinear form $\QQ$ on $\VV$, defined over $\R$, of parity $(-1)^k$,
whose associated flat Hermitian form $\QQ^h(\,\cdot\,,\,\cdot\,) :=
i^{-k}\, \QQ(\,\cdot\,,\,\bar{\cdot}\,)$ makes the
decomposition~\eqref{bundledecomposition} orthogonal and such that
$(-1)^p\QQ^h$ is positive definite on $\VV^{p,k-p}$.

Via parallel translation to a fixed fiber $V$ we may describe a
polarized variation of Hodge structure by a holomorphic \jdef{period
  map} $\Phi\colon M \to \DD/\Gamma$, where $\DD$ is the classifying
space of polarized Hodge structures on $V$ and $\Gamma$ is the
monodromy group.  We recall that $\DD$ is Zariski open in the smooth
projective variety \DC consisting of all filtrations $F$ in $V$, with
$\dim\,F^p = \sum_{r\geq p}\,h^{r,k-r}$, satisfying
$\,Q(F^p,F^{k-p+1}) = 0\, $, where $Q$ denotes the restriction of
$\QQ$ to $V$. The complex Lie group $\GC := \aut(V,Q)$ acts
transitively on $\DC$, and $\DD$ is an open orbit of $\GR :=
\aut(V_\R,Q)$.

Let $\jgg$ and $\lgr$ denote the Lie algebras of $\GC$ and $\GR$,
respectively.  The choice of a base point $F\in \DC$ defines a
filtration
\begin{equation*}
  F^a\jgg \ :=\ \{\,T\in\jgg\ :\ T\,F^p \subset F^{p+a}\,\}\,
\end{equation*}
compatible with the Lie bracket.  In particular, $F^{0}\jgg$ is the
isotropy subalgebra at $F$ and since $[F^{0}\jgg,F^{-1}\jgg]\subset
F^{-1}\jgg$, the quotient $F^{-1}\jgg/F^0\jgg$ defines a \GC-invariant
subbundle of the holomorphic tangent bundle of \DC ---the
\jdef{horizontal tangent bundle}.  Because of
\eqref{eq:horizontality}, the differential of $\Phi$ or, more
precisely, of any local lifting of $\Phi$ takes values on the
horizontal bundle.  Such maps are called \jdef{horizontal}.

Suppose now that $M$ has a smooth compactification $\conj{M}$ such
that $X := \conj{M}\setminus M$ is a normal crossings divisor. Around
a point of $X$, the \jdef{local variation} may be described by a
horizontal map
\begin{equation}\label{eq:localperiod}
  \Phi:(\Delta^*)^r \times \Delta^m\rightarrow \DD/\Gamma,
\end{equation}
where $\Delta$ is the unit disk in $\C$ and $\Delta^*$ the punctured
disk. We shall also denote by $\Phi$ its lifting to the universal
covering $U^r \times \Delta^m$, where $U$ is the upper-half plane.  We
let $z=(z_j)$, $t=(t_l)$ and $s=(s_j)$ be the coordinates on $U^r$,
$\Delta^m$ and $(\Delta^*)^r$ respectively. By definition, we have
$s_j = e^{2\pi i z_j}$.

Asymptotically, a period map has an algebraic component ---the
\jdef{nilpotent orbit}--- encoding the singularities of the connection
$\nabla$, and an analytic part described by a holomorphic map with
values in a nilpotent Lie algebra.  Assuming, for simplicity, that the
local monodromy of the variation is unipotent, let $N_1,\dots,N_r$
denote the monodromy logarithms.  Our convention is such that $\Phi(z
+ e_i, t) = \exp(N_i) \Phi(z , t)$, where $e_i$ denotes the $i$-th
standard vector.  It follows from Schmid's Nilpotent Orbit Theorem
\cite{ar:S-vhs} that the \DC-valued map
\begin{equation*}
  \Psi(s,t) \ :=\ \exp\biggl( -\sum_{j=1}^r \frac {\log s_j}{2\pi i}\ N_j
  \biggr)\cdot \Phi(s , t)
\end{equation*}
extends holomorphically to the origin.  The \jdef{limiting Hodge
  filtration} is $F_0 := \Psi(0,0)\in \DC$. The map
\begin{equation}
  \label{eq:nilpotent_orbit}
  \theta(z) \ :=\ \exp\biggl( \sum_{j=1}^r z_j\ N_j \biggr)\cdot F_0 \in \DC
\end{equation}
is holomorphic, horizontal, and \DD-valued for ${\rm Im}(z_j)\gg 0$;
\ie, is the period map of a local variation.

A nilpotent linear transformation $N\in\gl(V_\R)$ defines an
increasing filtration, the \jdef{weight filtration}, $W(N) $ of $V$,
defined over $\R$ and uniquely characterized by requiring that
$N(W_l(N))\subset W_{l-2}(N)$ and that $N^l:\gr_{l}^{W(N)}\rightarrow
\gr_{-l}^{W(N)}$ be an isomorphism.  It follows from
\cite[Theorem~3.3]{ar:CK-polarized} that if $N_1,\dots,N_r$ are the
monodromy logarithms of a local variation, then the weight filtration
$W(\sum \lambda_j N_j)$, $\lambda_j\in \R_{>0}$, is independent of the
choice of $\lambda_1,\dots,\lambda_r$ and, therefore, is associated
with the positive real cone $\CC \subset \jgg_\R$ spanned by
$N_1,\dots,N_r$.

The shifted weight filtration $W = W(\CC)[-k]$ and the limiting Hodge
filtration $F_0\in \DC$ define a \jdef{mixed Hodge Structure} on $V$;
\ie $F_0$ induces a Hodge structure of weight $\ell$ on $\gr_\ell^{W}$
for each $\ell$.  Recall (\cite[Theorem~2.13]{ar:CKS}) that mixed
Hodge structures are equivalent to (canonical) bigradings of $V$,
$I^{*,*}$, satisfying $I^{p,q}\equiv \conj{I^{q,p}}
\mod(\oplus_{a<p,b<q} I^{a,b})$. Thus, $W_l = \oplus_{p+q \leq l}
I^{p,q}$ and $F_0^a = \oplus_{p\geq a} I^{p,q}$.

A mixed Hodge structure $(W,F)$ is said to \jdef{split} over $\R$ if
$I^{p,q}= \conj{I^{q,p}}$; in that case the subspaces $V_l =
\oplus_{p+q = l} I^{p,q}$ define a real grading of $W$. A structure
for which $I^{p,q} = \{0\}$ if $p\neq q$ is said to be of
\jdef{Hodge-Tate} type.  A map $T \in \gl(V_\R)$ such that $T(I^{p,q})
\subset I^{p+a,q+b}$ is called a morphism of bidegree $(a,b)$.

A \jdef{polarized mixed Hodge structure}~\cite[(2.4)]{ar:CK-polarized}
of weight $k$ on $V_\R$ consists of a mixed Hodge structure $(W,F)$ on
$V$, a $(-1,-1)$ morphism $N\in \lgr$, and a nondegenerate,
$(-1)^k$-symmetric, bilinear form $Q$ such that
\begin{enumerate}
\item $N^{k+1}=0$,
\item $W = W(N)[-k]$, where $W[-k]_j = W_{j-k}$,
\item $Q(F^a,F^{k-a+1}) = 0$ and,
\item the Hodge structure of weight $k+l$ induced by $F$ on
  $\ker(N^{l+1}:\gr_{k+l}^{W}\rightarrow \gr_{k-l-2}^{W})$ is
  polarized by $Q(\cdot,N^l \cdot)$.
\end{enumerate}

It follows from Schmid's $SL_2$-orbit theorem \cite{ar:S-vhs} that the
mixed Hodge structure $(W(\CC)[-k],F_0)$ associated with a local
variation is polarized by every $N\in\CC$.  Conversely, given
commuting nilpotent elements $N_1,\ldots,N_r \in \lgr$ so that the
weight filtration $W(\sum \lambda_j N_j)$, $\lambda_j\in \R_{>0}$, is
independent of the choice of $\lambda_1,\dots,\lambda_r$, \ and $F_0
\in \DC$ such that $(W(\CC),F_0)$ is polarized by every element $N\in
\CC$, the map~\eqref{eq:nilpotent_orbit} is a period mapping for ${\rm
  Im}(z_j)$ sufficiently large \cite[Proposition 4.66]{ar:CKS}.
Moreover, if $(W(\CC),F_0)$ splits over $\R$, then $\theta(z) \in \DD$
for ${\rm Im}(z_j) > 0$.  We refer to the map $\theta$, or
equivalently, to $\{N_1,\dots,N_r;F_0\}$ as a \jdef{nilpotent orbit}.

The following example shows the relationship between nilpotent orbits
(equivalently, polarized mixed Hodge structures) and the Lefschetz
structure on the cohomology of a compact K\"ahler manifold.  This
point of view was introduced in \cite{ar:CKS2} where it was used to
obtain relations between the Lefschetz decompositions corresponding to
different K\"ahler classes.

\begin{example}
  \label{ex:MHS_for_kahler}
  If $X$ is a compact K\"ahler manifold of dimension $k$, the
  bigrading $I^{p,q}:= H^{k-q,k-p}(X)$ defines a mixed Hodge structure
  $(W,F)$ on $ H^*(X,\C)$ that splits over $\R$.  The interest of this
  construction lies in the fact that this mixed Hodge structure is
  polarized by the K\"ahler cone.  Indeed, the Hard Lefschetz Theorem
  is equivalent to the statement that if $\omega$ is a K\"ahler class
  and $L_\omega$ denotes multiplication by $\omega$, then $W =
  W(L_\omega)[-k]$; while the Hodge-Riemann bilinear relations imply
  that $L_\omega$ polarizes $(W,F)$ relative to the intersection form.
  The restriction of $(W,F)$ to $V := \oplus_{p=0}^k H^{p,p}$ defines
  a mixed Hodge structure of Hodge-Tate type.
\end{example}

We now describe the analytic component of a local variation.  The
bigrading associated with the limiting mixed Hodge structure $(W,F_0)$
defines a bigrading $I^{*,*}\jgg$ of the Lie algebra $\jgg$ by
$I^{a,b}\jgg :=\{X\in \jgg : X(I^{p,q})\subset I^{p+a,q+b} \}$.  Set
\begin{equation}
  \label{eq:def_pa}
  \jgp_a \ := \ \bigoplus_{q}I^{a,q}\jgg \quad \text{ and }\quad \jgg_-
  \ := \ \bigoplus_{a\leq
    -1}\jgp_a.
\end{equation}
The nilpotent subalgebra $\jgg_-$ is a complement of the stabilizer
subalgebra at $F_0$. Hence $(\jgg_-, X \mapsto \exp(X)\cdot F_0)$
provides a local model for the \GC-homogeneous space \DC near $F_0$.
Thus, locally around the origin, we may write $\Psi(s,t) =
\exp(\Gamma(s,t))\cdot F_0$, where $\Gamma(s,t)$ is a holomorphic
$\jgg_-$-valued map with $\Gamma(0,0) =0$.  We also write
\begin{equation*}
  \Phi(s,t) \ =\ \exp\biggl(\frac{1}{2\pi i} \sum_{j=1}^r \log(s_j)
  N_j\biggr)\cdot \exp(\Gamma(s,t))\cdot F_0 \ =\
  \exp\bigl(X(s,t)\bigr)\cdot F_0, 
\end{equation*}
where $X(s,t) \in \jgg_-$.  The horizontality of $\Phi$ now
translates, in terms of the gradings~\eqref{eq:def_pa}, into:
\begin{equation}
  \label{eq:horiz}
  \exp\bigl(-X(s,t)\bigr)\, d\exp\bigl(X(s,t)\bigr)\ =\ dX_{-1} \in
  \jgp_{-1}\otimes T^*((\Delta^*)^r\times\Delta^m),
\end{equation}
where $X_{-1}$ denotes the $\jgp_{-1}$-graded part of $X$.  In
particular,
\begin{equation} 
  \label{eq:integcond}
  dX_{-1}\wedge dX_{-1}\ =\ 0,
\end{equation}
where $X_{-1} = \frac{1}{2\pi i} \sum_{j=1}^r \log(s_j) N_j +
\Gamma_{-1}$.

The following result, which follows from
\cite[Theorem~2.8]{ar:CK-luminy} and
\cite[Theorem~2.7]{ar:CF-asymptotics}, shows that the nilpotent orbit
together with the $\jgp_{-1}$-valued holomorphic function
$\Gamma_{-1}$ completely determine the local variation:

\begin{theorem}\label{th:improved_2.8}
  Let $\{N_1,\ldots,N_r;F_0\}$ be a nilpotent orbit and $R:\Delta^r
  \times \Delta^m \rightarrow \jgp_{-1}$ be a holomorphic map with
  $R(0,0)=0$.  Define $X_{-1}(z,t) := \sum_{j=1}^r z_j N_j+R(s,t)$,
  $s_j = e^{2\pi i z_j}$, and suppose that the differential equation
  \eqref{eq:integcond} holds.  Then, there exists a unique period
  mapping
  \begin{equation*}
    \Phi(s,t) \ =\ \exp\biggl(\frac{1}{2\pi i} \sum_{j=1}^r \log(s_j)
    N_j\biggr)\cdot \exp(\Gamma(s,t))\cdot F_0,
  \end{equation*}
  defined in a neighborhood of the origin in $\Delta^{r+m}$ such that
  $\Gamma_{-1} = R$.
\end{theorem}

In the ensuing sections we will be concerned with a special type of
maximally degenerating variation. These are relevant to the study of
mirror symmetry and, from a Hodge theoretic perspective they have the
advantage of allowing us to use a canonical system of coordinates on
the parameter space of the variation. Following Morrison~\cite[Def.
3]{ar:morr-compactifications}, we consider

\begin{definition}\label{maxunip}
  Given a polarized variation of Hodge structure of weight $k$ over
  $(\Delta^*)^r$ whose monodromy is unipotent, we say that $0\in
  \Delta^r$ is a \jdef{maximally unipotent boundary point} if
  \begin{enumerate}
  \item $\dim I^{k,k} = 1$, $\dim I^{k-1,k-1} = r$ and $\dim I^{k,k-1}
    = \dim I^{k-2,k}= 0$, where $I^{*,*}$ is the bigrading associated
    to the limiting mixed Hodge structure and,
  \item $\vspan_\C\{N_1(I^{k,k}),\ldots,N_r(I^{k,k})\} = I^{k-1,k-1}$,
    where $N_j$ are the monodromy logarithms of the variation.
  \end{enumerate}
\end{definition}

The limiting Hodge filtration $F_0$ and the holomorphic function
$\Gamma$ of a local variation depend on the choice of coordinates on
$(\Delta^*)^r$.  However, in the maximally unipotent case we may
normalize our choices as follows.

\begin{proposition}
  \label{prop:canonical_coordinates}
  Let $\Phi = \exp(\sum_{j=1}^r \frac{1}{2\pi i}\log(s_j) N_j)\cdot
  \exp(\Gamma(s))\cdot F_0$ be a polarized variation of Hodge
  structure that has a maximally unipotent boundary point at $0\in
  \Delta^r$. Then, there is a coordinate system on $\Delta^r$, unique
  up to scaling, where $\Gamma$ satisfies $\Gamma(I^{1,1}) = 0$.
\end{proposition}

For a proof of Proposition~\ref{prop:canonical_coordinates},
see~\cite[\S 3]{ar:CF-asymptotics}.  We will refer to these as
\jdef{canonical coordinates}. They are standard in the physics
literature and their Hodge-theoretic interpretation is due to
D.~Morrison~\cite{ar:morr-picard-fuchs} and
P.~Deligne~\cite{ar:Del-local_behavior}.



\section{Frobenius modules}
\label{sec:frobenius_modules}

The cohomology of even degree of a compact manifold is a graded
Frobenius algebra relative to cup product and the intersection form.
When $X$ is K\"ahler, the Hard Lefschetz Theorem and the Hodge-Riemann
bilinear relations describe the action of $H^{1,1}(X)$ on the full
cohomology.  We abstract these properties in the notion of a (framed)
Frobenius module.

Let $V=\oplus_{p=0}^k V_{2p}$ be a graded \C-vector space and \CB a
symmetric nondegenerate bilinear form on $V$ pairing $V_{2p}$ with
$V_{2(k-p)}$. Let $\{T_a\}_{0\leq a\leq m}$ be a \CB-self dual, graded
basis of $V$.  We will refer to $\{T_a\}$ as an \jdef{adapted basis}.
For $0\leq a\leq m$ define $\delta(a)$ by $\CB(T_{\delta(a)},T_b) =
\delta_{a b}$ for all $b=0,\ldots,m$.  We also set $\ti{a} := p$ if
and only if $ T_a \in V_p$ and assume that the map $\sim\ :
\{0,\ldots,m\} \rightarrow \{0,\ldots, 2k\}$ is increasing.

\begin{definition}
  $(V,\CB,e,*)$ is a \jdef{graded $V_2$-Frobenius module} of weight
  $k$ if
  \begin{enumerate}
  \item $e\neq 0$ and $V_0 = \langle e \rangle$.
  \item $V$ is a graded $\sym V_2$-module under $*$.
  \item For all $v_1, v_2 \in V$ and $w \in V_2$
    \begin{equation}
      \label{eq:frobenius_condition}
      \CB(w * v_1, v_2) \ =\ \CB(v_1,w * v_2)
    \end{equation}
  \item $w * e = w$ for all $w\in V_2$.
  \end{enumerate}
\end{definition}

Since $T_0\in V_0$, it must be a non-zero multiple of $e$ and we
assume that an adapted basis satisfies $T_0 = e$.  Clearly, the fact
that $V$ is a $\sym V_2$-module is equivalent to
\begin{equation}
  \label{eq:sym_condition}
  T_j * (T_l * T) \ =\ T_l * (T_j * T) \text{ for all }
  T_j, T_l \in V_2 \text{ and } T\in V.
\end{equation}

We say that $V$ is \jdef{real} if $V$ has a real structure, $V_\R$,
compatible with its grading, $*$ is real, $e\in V_\R$, and $\CB$ is
defined over $\R$.

\begin{example}
  \label{ex:FM_for_kahler}
  If $X$ is a compact K\"ahler manifold of dimension $k$, let
  $V_{2p}:= H^{p,p}(X)$, $\CB_{int}$ the intersection pairing on
  $V:=\oplus_{p=0}^k V_{2p}$, and $\smallsmile$ the restriction of the cup
  product to $V$. Then, $(V,\CB_{int},1,\smallsmile)$ defines a real
  Frobenius module. The real structure is induced by $H^*(X,\R)$.
\end{example}

As in the case of the cohomology of a compact K\"ahler manifold, to
any real Frobenius module we can associate a Hodge-Tate mixed Hodge
structure:
\begin{equation}
  \label{eq:bigrading_from_grading}
  I^{p,p} \ :=\  V_{2(k-p)}.
\end{equation}

The multiplication operator $L_w\in\jend(V)$, $w\in V_2$, is an
infinitesimal automorphism of the bilinear form
\begin{equation}
  \label{eq:Q_from_B_again}
  Q(v_a,v_b) \ :=\ (-1)^{k+\ti{a}/2} \CB(v_a,v_b), 
\end{equation}
as well as a $(-1,-1)$-morphism of the associated mixed Hodge
structure.  We will say that $w\in V_2\cap V_\R$ \jdef{polarizes} $V$
if the mixed Hodge structure $(I^{*,*},Q,L_w)$ is polarized.  A real
Frobenius module $V$ is said to be \jdef{polarizable} if it contains a
polarizing element. Given a polarizing element $w$, the set of
polarizing elements is an open cone in $V_{2} \cap V_\R$. We can then
choose a basis $T_1, \ldots, T_r$ of $V_{2}\cap V_\R$ spanning a
simplicial cone $\CC$ contained in the closure of the polarizing cone
and with $w\in \CC$.  Such a choice of a basis of $V_2$ will be called
a \jdef{framing} of the polarized Frobenius module.

Given an adapted basis $\{T_0,\ldots,T_m\}$ of $V$, let
$z_0,\ldots,z_m$ be the corresponding linear coordinates on $V$ and
set $q_j:= \exp(2\pi i z_j)$ for $j=1,\ldots, r:=\dim V_2$ .  We may
identify $U^r \cong (V_2 \cap V_\R) \oplus i\,\CC$ and view
the correspondence
\begin{equation*}
  \sum_{j=1}^r z_j T_j \in (V_2 \cap V_\R) \oplus i\,\CC
  \mapsto (q_1,\dots,q_r) \in (\Delta^*)^r
\end{equation*}
as the natural covering map.

\begin{proposition}
  \label{prop:framed_FM_NO}
  Framed, real Frobenius modules of weight $k$ are equivalent to
  nilpotent orbits of weight $k$ whose limiting mixed Hodge structure
  is of Hodge-Tate type, split over $\R$, have a marked real element
  in $F^k$, and have the origin as a maximally unipotent boundary
  point.
\end{proposition}

\begin{proof}  
  Let $(V,\CB,e,*)$ be a real Frobenius module with framing
  $T_1,\dots,T_r$.  Set $N_j := L_{T_j}$ and $F^p := \oplus_{a\geq p}
  I^{a,a}$.  Then $\{N_1,\dots,N_r;F \}$ is a nilpotent orbit.  The
  element $e \in I^{k,k} = F^k$ is a distinguished real element and
  the conditions of Definition~\ref{maxunip} are clearly satisfied.
  
  Conversely, suppose $\{N_1,\dots,N_r;F \}$ is a nilpotent orbit
  whose limiting mixed Hodge structure is of Hodge-Tate type, split
  over $\R$ and satisfies both conditions of Definition~\ref{maxunip}.
  Set $V_{2p}:= I^{k-p,k-p}$; in particular, the marked element $e \in
  F^k = I^{k,k} = V_0$ and it follows from (2) in
  Definition~\ref{maxunip} that the map 
  \begin{equation*}
    N \in \vspan_\C\{N_1,\dots,N_r\} \mapsto N(e)
  \end{equation*}
  identifies the polynomial algebra $\C[N_1,\dots,N_r]$ with $\sym
  V_2$ and defines a $\sym V_2$-action on $V$. Let \CB be defined from
  the polarization $Q$ as in \eqref{eq:Q_from_B_again}, then since the
  monodromy transformations $N_j$ are infinitesimal automorphisms of
  $Q$, \eqref{eq:frobenius_condition} is satisfied.  Thus,
  $(V,\CB,e,*)$ is a Frobenius module.  The equivalence between
  nilpotent orbits and polarized mixed Hodge structures implies that
  $T_j = N_j(e)$, $j=1,\dots,r$, are a framing of $V$ and the fact
  that $N_1,\dots,N_r$ are real implies that the Frobenius structure
  is real.
\end{proof}

A Frobenius module structure may also be encoded in a polynomial of
degree $3$ in the variables $z_0,\dots,z_m$.  Indeed, if we let
\begin{equation*}
  \phi_0(z_0,\ldots,z_m) \ :=\ \sum_{\ti{j}=2,\ 0\leq\ti{a},\ti{b}\leq
  2k} z_j z_a z_b\,C(\ti{a})\,\CB(T_j * T_a, T_b)\,, 
\end{equation*}
with
\begin{equation*}
  C(\ti{a})\ :=\
  \begin{cases}
    \frac{1}{6} \text{ if } k=3 \text{ and } \ti{a}=2,\\
    \frac{1}{4} \text{ if } k\neq 3 \text{ and } \ti{a}=2 \text{ or }
    \ti{a}=2k-4,\\
    \frac{1}{2} \text{ otherwise},
  \end{cases}
\end{equation*}
then we recover the $\sym V_2$-action by:
\begin{equation*}
  T_j * T_a \ :=\ \sum_{\ti{c}=\ti{a}+2}
  \pdiii{\phi_0}{j}{a}{\delta(c)} T_c\ ;\quad j=1,\dots,r\,.
\end{equation*}
The polynomial $\phi_0$ is called a \jdef{(classical) potential} for
the Frobenius module.

We may generalize this construction by considering deformations of the
classical potential.  This is motivated by the construction of the
quantum product as a deformation of the cup product on the cohomology.
We assume, for simplicity, that $k>3$.  Let
$R:=\C\{q_1,\ldots,q_r\}_0$ denote the ring of convergent power series
vanishing for $q_1 = \cdots = q_r = 0$ and $R'$ be its image under the
map induced by $q_j\mapsto e^{2\pi i z_j}$ for $1\leq j \leq r$.

\begin{definition}
  \label{def:quantum_potential}
  Let $(V,\CB,e,*)$ be a Frobenius module of weight $k > 3$ with
  classical potential $\phi_0$. A \jdef{quantum potential} on $V$ is a
  function $\phi:V\rightarrow \C$ of the form $\phi = \phi_0 +
  \phi_\h$, where
  \begin{equation}
    \label{eq:quantum_potential}
    \phi_\h(z) \ :=\ \sum_{\ti{a}=2k-4} z_a
    \phi_h^a(z_1,\ldots,z_r) +\sum_{\substack{2 < \ti{a}< 2k-4\\
    \ti{a}+\ti{b} = 2k-2}} z_a z_b \phi_h^{a b}(z_1,\ldots,z_r), 
  \end{equation}
  with $\phi_h^a, \phi_h^{a b}\in R'$ and such that
  \begin{equation}
    \label{eq:WDVV_graded}
    \sum_{\ti{c}=\ti{a}+2} \pdiii{\phi}{l}{a}{\delta(c)}\
    \pdiii{\phi}{j}{c}{\delta(d)} \ =\ \sum_{\ti{c}=\ti{a}+2}
    \pdiii{\phi}{j}{a}{\delta(c)} \ \pdiii{\phi}{l}{c}{\delta(d)}
  \end{equation}
  holds for all $a$, $\ti{j} = \ti{l} = 2$ and $\ti{d} = \ti{a}+4$.
\end{definition}

Given a quantum potential $\phi$ on $(V,\CB,e,*)$, we can define a
deformation of the module structure by
\begin{equation}
  \label{eq:product_from_potential}
  T_j \cdot_q T_a \ :=\ \sum_{\ti{c}=\ti{a}+2}
  \pdiii{\phi}{j}{a}{\delta(c)} T_c\,,\quad \text{ with }
  q=(q_1,\ldots,q_r) \in \Delta^r\,.
\end{equation}
We should stress that, even though the right side
of~\eqref{eq:product_from_potential} depends explicitly on the
variables $z_0,\ldots, z_m$, \eqref{eq:quantum_potential} implies that
it is actually a function of $q_1,\ldots,q_r$.
Condition~\eqref{eq:WDVV_graded} guarantees
that~\eqref{eq:product_from_potential} defines an action of $\sym V_2$
for all $q$. Moreover, $(V,\CB,T_0,\cdot_q)$ is a Frobenius module of
weight $k$ for all $q$, and $\cdot_0 = *$. We will say that a
deformation of the Frobenius module $V$ is \jdef{framed} if $V$ is
framed.

\begin{remark} 
  Definition~\ref{def:quantum_potential} abstracts the properties of
  the graded portion of the Gromov-Witten potential needed to describe
  the action of $H^{1,1}(X,\C)$ in the small quantum cohomology ring
  of a Calabi-Yau manifold $X$.  In particular, \eqref{eq:WDVV_graded}
  is a graded component of the WDVV equations.  We refer to ~\cite[\S
  8.2, \S 8.3]{bo:CK-mirror} and~\cite[\S 5]{ar:CF-asymptotics} for
  details.
  
  We can extend the definition of quantum potential to the weight $3$
  case by taking $\phi = \phi_0 + \phi_\h$ for $\phi_\h \in R'$. With
  this notion, all the results from Sections~\ref{sec:correspondence}
  and~\ref{sec:a_model_variation} extend to this weight. For $V$ of
  weight $1$ or $2$, the Frobenius module is determined by \CB and
  $e$; hence no deformations are possible.
\end{remark}



\section{Correspondence}
\label{sec:correspondence}

In this section we will prove the main result of this paper, namely
the correspondence between deformations of framed Frobenius modules
and degenerating polarized variations of Hodge structures. In
\S\ref{sec:a_model_variation} we will show that when the deformation
arises from the quantum product of a Calabi-Yau manifold, the
associated variation of Hodge structure is the so-called $A$-model
variation.

\begin{theorem}
  \label{th:correspondence_grl}
  There is a one to one correspondence between
  \begin{itemize}
  \item Quantum potentials $\phi$ on a framed Frobenius module
    $(V,\CB,e,*)$ of weight $k$, and
  \item Germs of polarized variations of Hodge structure of weight $k$
    on $V$ degenerating at a maximally unipotent boundary point to a
    limiting mixed Hodge structure of Hodge-Tate type, split over
    $\R$, and together with a marked real point $e\in F^k$.
  \end{itemize}
  
  Under this correspondence, classical potentials ---equivalently,
  framed Frobenius modules--- correspond to nilpotent orbits as in
  Proposition~\ref{prop:framed_FM_NO}.
\end{theorem}

\begin{proof}
  Let $(V,\CB,e,*)$ be a framed Frobenius module of weight $k$,
  $\{T_0,\ldots,T_m\}$ an adapted basis, and let
  $\{N_1,\ldots,N_r;F\}$ be the nilpotent orbit associated by
  Proposition~\ref{prop:framed_FM_NO}.  Given a quantum potential
  $\phi = \phi_0 + \phi_\h$ on $V$ define
  \begin{equation}
    \label{eq:gamma_from_potential}
    \Gamma_{-1}(q)(T_a) \ :=\ \sum_{\ti{c}=\ti{a}+2}
    \pdii{\phi_\h(q)}{a}{\delta(c)} T_c.
  \end{equation}
  Notice that because of \eqref{eq:quantum_potential}, $\Gamma_{-1}$
  is holomorphic on some open neighborhood of $q=0\in\Delta^r$,
  $\Gamma_{-1}(0)=0$, and it takes values on $\jgp_{-1}$ relative to
  the grading~\eqref{eq:def_pa} defined by the limiting mixed Hodge
  structure of $\{N_1,\ldots,N_r;F\}$.
  
  As before, we set $X_{-1}(q) := \frac{1}{2\pi i}\sum_{j=1}^r
  \log(q_j) N_j +\Gamma_{-1}(q) \in \jgp_{-1}$ and note that the
  deformed Frobenius structure may be recovered from $X_{-1}(q)$ by
  \begin{equation}
    \label{eq:product_from_X}
    T_j\cdot_q T_a\ =\ \pd{X_{-1}}{z_j}(T_a)\,; \quad \ti{j}=2\,,\
    0\leq a\leq m\,.
  \end{equation}
  
  The equations~\eqref{eq:WDVV_graded} imply that $X_{-1}$ satisfies
  the integrability condition~\eqref{eq:integcond}.  Indeed,
  \begin{equation}
    \label{eq:horizontality_sym}
    \begin{split}
      dX_{-1} \wedge dX_{-1} \ =\ 0 & \iff \pd{X_{-1}}{z_j}
      \pd{X_{-1}}{z_l}\ =\ \pd{X_{-1}}{z_l} \pd{X_{-1}}{z_j} \\ & \iff
      T_j\cdot_q(T_l\cdot_q T_a)\ =\ T_l\cdot_q(T_j\cdot_q T_a)\,,
    \end{split}
  \end{equation}
  which, by~\eqref{eq:WDVV_graded}, holds whenever $\ti{j} = \ti{l} =
  2$ and all $a$.  Theorem~\ref{th:improved_2.8} now implies that
  $X_{-1}$ defines a unique polarized variation of Hodge structure on
  a neighborhood of $0\in \Delta^r$ whose nilpotent orbit is
  $\{N_1,\ldots,N_r;F\}$.  Hence the origin is a maximally unipotent
  boundary point and the limiting mixed Hodge structure is of
  Hodge-Tate type.
  
  Conversely, let $\Phi$ be the period map of a local variation having
  a maximally unipotent boundary point at the origin.  Let
  $\{N_1,\ldots,N_r;F\}$ be the corresponding nilpotent orbit and
  $I^{*,*}$ the limiting mixed Hodge structure, which we assume to be
  of Hodge-Tate type.  Let $(V,\CB,e,*)$ be the real, framed Frobenius
  module given by Proposition~\ref{prop:framed_FM_NO} and $\phi_0$ the
  corresponding classical potential.  Let $\{T_0,\ldots,T_m\}$ be an
  adapted basis such that $T_j = N_j(e)$, $j=1,\dots,r$. Using
  canonical coordinates $q$ on $\Delta^r$,  we define a
  quantum potential from the holomorphic function $\ \Gamma \colon
  \Delta^r \to \jgg_{-}\,$ associated with $\Phi$ by:
  \begin{eqnarray*}
    \phi_h^{a b}(q) &:=& \frac{1}{2}\, \CB(\Gamma_{-1}(T_a),T_b)
    \text{ for } 2<\ti{a}< 2k-4 \text{ and } \ti{a}+\ti{b}=2k-2\\  
    \phi_h^a(q) &:=& \CB(-\Gamma_{-2}(T_a),T_0) \text{ for }
    \ti{a}=2k-4\\
    \phi_\h &:=& \sum_{\ti{a}\ =\ 2k-4} z_a \phi_h^a +
    \sum_{\substack{2\ <\ \ti{a}\ <\ 2k-4\\ \ti{a}+\ti{b}\ =\ 2k-2}}
    z_a z_b \phi_h^{a b}\\
    \phi  &:=& \phi_0 + \phi_\h. 
  \end{eqnarray*}
  Clearly, $\phi_\h$ is as in \eqref{eq:quantum_potential}.  In order
  to verify that~\eqref{eq:WDVV_graded} is satisfied we consider the
  associated deformation~\eqref{eq:product_from_potential} of the
  Frobenius module structure
  \begin{equation*}
    T_j \cdot_q T_a \ :=\ \sum_{\ti{c}=\ti{a}+2}
    \pdiii{\phi}{j}{a}{\delta(c)} T_c
  \end{equation*}
  and show that it may also be given as 
  \begin{equation}\label{equality}
    T_j \cdot_q T_a = \frac {\partial X_{-1}}{\partial z_j}(T_a).
  \end{equation}
  Indeed, for $2<\ti{a}<2k-4$ we have $\Gamma_{-1}(T_a) =
  \sum_{\ti{c}=\ti{a}+2} \phi_h^{a\delta(c)} T_c$, so that
  \begin{equation*}
    \begin{split}
      \pd{\Gamma_{-1}}{z_j}(T_a) \ &=\ \sum_{\ti{c}=\ti{a}+2}
      \pd{}{z_j}\phi_h^{a\delta(c)} T_c \ =\ \sum_{\ti{c}=\ti{a}+2}
      \pdiii{}{j}{a}{\delta(c)} \sum_{\ti{u}+\ti{v}=2k-2} \frac{1}{2}
      \,z_u z_v\, \phi_h^{u v}\\
      &=\ \sum_{\ti{c}=\ti{a}+2}
      \pdiii{ \phi_\h}{j}{a}{\delta(c)} T_c,
    \end{split}
  \end{equation*}
  where we have used that $\phi_h^{a b} = \phi_h^{b a}$. Then
  \begin{equation*}
    \begin{split}
      \pd{X_{-1}}{z_j}(T_a) &=\  N_j(T_a)\  +\ 
      \pd{\Gamma_{-1}}{z_j}(T_a)\\ &=\ \sum_{\ti{c} = \ti{a}+2}
      \pdiii{\phi_0}{j}{a}{\delta(c)} T_c \ +\  \sum_{\ti{c}=\ti{a}+2}
      \pdiii{ \phi_\h}{j}{a}{\delta(c)} T_c\\
      &=\ \sum_{\ti{c}=\ti{a}+2}
      \pdiii{ \phi}{j}{a}{\delta(c)} T_c \ =\ T_j\cdot_q T_a.
    \end{split}
  \end{equation*}
  In order to verify~\eqref{equality} when $\ti{a}=2k-4$ we first
  prove the identity
  \begin{equation}
    \label{eq:original_formula}
    \Gamma_{-1}(T_a)\ =\ \sum_{\ti{c}=2k-2} \pd{}{z_{\delta(c)}}
    \CB(-\Gamma_{-2}(T_a),T_0) \ T_c\,,\ \ \ti{a} = 2k-4
  \end{equation}
  as a consequence of the horizontality condition~\eqref{eq:horiz}.
  If this condition is rewritten in terms of $G(q) :=\exp \Gamma(q)$
  and $\Theta = \sum N_j\, dz_j$ we get
  \begin{equation*}
    dG \ =\ [G,\Theta] + G\, d\Gamma_{-1}.
  \end{equation*}
  This equation is graded by~\eqref{eq:def_pa} and its homogeneous
  pieces are
  \begin{equation}
    \label{eq:horiz_graded}
    dG_{-\ell} \ =\ [G_{-\ell+1},\Theta] +G_{-\ell+1}\,
    d\Gamma_{-1},\quad \ell\geq 2. 
  \end{equation}
  In particular, for $\ell=2$ we obtain
  \begin{equation*}
    d\Gamma_{-2} \ =\ [\Gamma_{-1},\Theta+\frac{1}{2}\, d\Gamma_{-1}].
  \end{equation*}
  Evaluating at $T_a$ and given that the canonical coordinates
  $(q_1,\dots,q_r)$ are characterized by $\Gamma_{-1}(T_b)=0$ for all
  $\ti{b}=2k-2$, we obtain
  \begin{equation}
    \label{eq:horiz_2}
    d\Gamma_{-2}(T_a) \ =\ -\Theta \bigl(\Gamma_{-1}(T_a)\bigr).
  \end{equation}
  By the \CB-self-duality of the basis $\{T_0,\dots,T_m\}$, we can
  write
  \begin{equation}
    \label{eq:Gamma_expansion}
    \Gamma_{-1}(T_a) \ =\ \sum_{\ti{c}=2k-2}
    \CB(\Gamma_{-1}(T_a),T_{\delta(c)}) T_c.
  \end{equation}
  Now, if $\ti{c}=2k-2$ and $j=1,\dots,r$, then $N_j(T_c)=\delta_{j
    c}T_m$ and, therefore, $\Theta(T_c) = T_m\, dz_{\delta(c)}$
  and~\eqref{eq:horiz_2},~\eqref{eq:Gamma_expansion} imply
  \begin{equation*}
    d\Gamma_{-2}(T_a) \ =\ -\sum_{\ti{c}=2k-2} dz_{\delta(c)}
    \CB(\Gamma_{-1}(T_a),T_{\delta(c)}) T_m,
  \end{equation*}
  so that,
  \begin{equation*}
    \pd{}{z_{\delta(c)}} \Gamma_{-2}(T_a) \ =\
    -\CB(\Gamma_{-1}(T_a),T_{\delta(c)}) T_m 
  \end{equation*}
  implying that
  \begin{equation}
    \label{eq:final_relation}
    \CB\biggl(\pd{}{z_{\delta(c)}} \Gamma_{-2}(T_a),T_0\biggr) \ =\
    -\CB(\Gamma_{-1}(T_a),T_{\delta(c)}). 
  \end{equation}
  Finally,~\eqref{eq:original_formula} follows from
  applying~\eqref{eq:final_relation} to~\eqref{eq:Gamma_expansion}.

  Thus, if $\ti{a} = 2k-4$,
  \begin{equation*}
    \begin{split}
      \pd{\Gamma_{-1}}{z_j}(T_a) \ &=\ \sum_{\ti{c}=2k-2}
      \pd{}{z_{j}}\,\pd{}{z_{\delta(c)}}\,
      \CB(-\Gamma_{-2}(T_a),T_0) \ T_c\\
      &=\ \sum_{\ti{c}=\ti{a} + 2}
      \pd{}{z_{j}}\,\pd{}{z_{\delta(c)}}\,\pd{}{z_{a}}\,
      \biggl(\sum_{\ti{b} = 2k-4} z_b \phi_h^b(q)\biggr) \ T_c\\
      &=\ \sum_{\ti{c}=\ti{a} + 2} \pdiii{\phi_\h}{j}{a}{\delta(c)}
      T_c\,.
    \end{split}
  \end{equation*}
  and \eqref{equality} follows as before.
  
  Given \eqref{equality}, the equivalences in
  \eqref{eq:horizontality_sym} show that the integrability
  condition~\eqref{eq:integcond} implies that the quantum potential
  $\phi$ satisfies~\eqref{eq:WDVV_graded}.
  
  Finally, we note that \eqref{equality} and \eqref{eq:product_from_X}
  imply that these correspondences are inverses of each other.
\end{proof}



\section{$A$-model variation}
\label{sec:a_model_variation}

Here we will show that the polarized variation of Hodge structure
associated to a quantum potential by
Theorem~\ref{th:correspondence_grl} agrees with the $A$-model
variation defined, in the case of the cohomology on a Calabi-Yau
manifold, by the Gromov-Witten potential, as in, for example,
\cite[Chapter 8]{bo:CK-mirror}.  As a byproduct we give a different
proof of the fact that the $A$-model variation associated with a
general potential, in the sense of
Definition~\ref{def:quantum_potential}, is a polarized variation of
Hodge structure.

We begin by recalling the definition of the $A$-model variation.  Let
$\phi = \phi_0 + \phi_\h$ be a quantum potential on the framed
Frobenius module $(V,\CB,e,*)$. Let $\{T_0,\ldots,T_m\}$ be an adapted
basis of $V$ and $(z_0,\dots,z_m)$ the corresponding linear
coordinates on $V$; set $q_j=\exp(2\pi i z_j)$ for $j=1,\dots,r$.  We
view $(q_1,\dots,q_r)$ as coordinates on $(\Delta^*)^r$.  Let $\nabla$
be the connection on the vector bundle $\VV := (\Delta^*)^r
\times V$ defined on a constant section $T$ by
\begin{equation}
  \label{eq:nabla_definition}
  \nabla_{\pd{}{q_j}} T\ :=\ \frac{1}{2\pi i q_j} T_j \cdot_q T.
\end{equation}

\begin{proposition}
  The connection $\nabla$ is flat.  It has a simple pole at $q_j=0$
  and its residue is the nilpotent operator
  \begin{equation}
    \label{eq:residue}
    \res_{q_j=0}(\nabla)(T_a)\ =\
    \frac{1}{2\pi i} \left(\sum_{\ti{c}=\ti{a}+2}
      \pdiii{\phi_0}{j}{a}{\delta(c)} T_c\right). 
  \end{equation}
\end{proposition}

\begin{proof}
  Given the definition of the quantum product
  \eqref{eq:product_from_potential} and \eqref{eq:nabla_definition},
  if $T_a$ denotes a constant section,
  \begin{equation*}
    \nabla_{\pd{}{q_j}} T_a\ =\ \frac{1}{2\pi i q_j}
    \left(\sum_{\ti{c}=\ti{a}+2} \pdiii{\phi_0}{j}{a}{\delta(c)}
    T_c\right) + H_{j a}(q)
  \end{equation*}
  for some function $H$, which extends holomorphically to $0\in
  \Delta^r$.  This implies the residue assertion.
  
  The curvature of $\nabla$ reduces to
  \begin{equation*}
    R_\nabla\left(\pd{}{q_j},\pd{}{q_l}\right) (T_a) = \frac{1}{2\pi i}
    \left(\frac{1}{q_l} \nabla_{\pd{}{q_j}}(T_l\cdot T_a) -\frac{1}{q_j}
    \nabla_{\pd{}{q_l}}(T_j\cdot T_a)\right).
  \end{equation*}
  A straightforward computation shows that this last expression
  vanishes since $\phi$ satisfies~\eqref{eq:WDVV_graded}.
\end{proof}

\begin{remark}
  It follows from~\eqref{eq:residue} that the operators
  $\res_{q_j=0}(\nabla)$ agree, up to a constant, with the morphisms
  $L_{T_j}$ of left multiplication by $T_j$ in the Frobenius module
  $(V,\CB,e,*)$.
\end{remark}

Consider the flags of subbundles of $\VV$:
\begin{equation*}
  \FF^p\ :=\ (\Delta^*)^r\times (\oplus_{a\geq p}
  V_{2(k-a)}) \quad \text{ and }\quad \UU_\ell\ :=\
  (\Delta^*)^r \times (\oplus_{b\geq\ell} V_{2b}).  
\end{equation*}

\begin{proposition}
  \label{prop:sharpsections}
  The subbundles $\FF^p$ satisfy Griffiths'
  horizontality~\eqref{eq:horizontality}.  Moreover, for any given
  $\hat{q}\in(\Delta^*)^r$, there is a (multivalued) flat frame of
  $\VV$, $\{T_a^\flat\}$, such that $T_a^\flat(q) \equiv T_a
  \mod{\UU_{\ti{a}+1}}$ and $T_a^\flat(\hat{q}) = T_a$.
\end{proposition}

\begin{proof} 
  Since the maps $\ T \mapsto T_j \cdot_q T\ $ are homogeneous of
  degree $2$, the horizontality follows directly
  from~\eqref{eq:nabla_definition}.
  
  Since $\nabla$ defines a connection on the bundle $\UU_\ell$
  inducing a trivial connection on $\UU_\ell/\UU_{\ell+1}$, the second
  statement follows.
\end{proof}

Next, we want to compute the monodromy of $\nabla$. We fix all the
coordinates $q_i$ for $i\neq j$ and consider the one-dimensional
problem around $q_j=0$. The flat sections $T_a^\flat$ can be written
in terms of the constant sections as $T_a^\flat = \sum_b f_{b a} T_b$,
and the flatness condition leads to the ODE with a regular singularity
at the origin
\begin{equation}
  \label{eq:flatness_coordinates}
  \pd{f_{ba}}{q_j}\ =\ - \sum_c \left( \frac{1}{q_j}
    (\res_{q_j=0}(\nabla))_{b c} + H_{j c \delta(b)}\right) f_{ca},
\end{equation}
where $H_{j c d}$ are holomorphic at $q_j=0$. Therefore, classical
results for such an equation (see~\cite[Ch. 4, Thm.  4.1]{bo:CL-ODE})
imply that the coefficients $f_{b a}$ are of the form
\begin{equation}
  \label{eq:flat_coordinates}
  f_{ba}(q)\ =\ \left(G(q_j)\, \exp(- \log(q_j)
  \res_{q_j=0}(\nabla))\right)_{ba} 
\end{equation}
for some function $G$, holomorphic at $q_j=0$, with $G(0) = \idM_n$ .

Parallel transport around $q_j=0$, in the anti-clockwise direction,
gives that the monodromy of $\nabla$, written relative to the frame
$\{T_a^\flat\}$, is 
\begin{equation*} 
  M_j\ :=\ \exp\bigl(- 2\pi i\res_{q_j=0}(\nabla)\bigr).
\end{equation*}
We let $N_j:= -\log(M_j) = 2\pi i\res_{q_j=0}(\nabla)$.  Notice that,
in view of~\eqref{eq:residue}, the monodromy in a flat frame can be
computed purely in terms of the classical potential. All together we
conclude:

\begin{proposition}
  \label{prop:N_product_0}
  The matrix of the local monodromy logarithm operator $N_j$ with
  respect to the frame $\{T_a^\flat\}$ coincides with the matrix of
  left $*$-multiplication by $T_j$, $L_{T_j}$, with
  respect to the basis $\{T_a\}$.
\end{proposition}

The fact that $\nabla$ has a simple pole at $q_j=0$ with nilpotent
residue $L_{T_j}$ allows us to construct Deligne's canonical extension
$(\VV^c, \nabla^c)$ \cite{bo:Del-equations} which is characterized by
the fact that
\begin{equation}
  \label{eq:flatc_sections}
  \ti{T}_a\ :=\ \exp\biggl(\sum_{j=1}^r\frac{\log(q_j)}{2\pi i} N_j\biggr)
  T_a^\flat\,,\quad a = 0,\dots,m, 
\end{equation}
are a flat frame of $(\VV^c, \nabla^c)$.

\begin{proposition}
  \label{prop:david_854}  
  For $a=0,\ldots,m$, $\ti{T}_a$ is the unique $\nabla^c$-flat section
  of $\VV^c$ such that $\ti{T}_a \equiv T_a \mod \UU_{\ti{a}+1}$, and
  $\ti{T}_a(\hat q) = T_a$.  The matrix of $N_j$ acting on the frame
  $\{\ti{T}_a\}$ equals the matrix of the classical product $*$ acting
  on $\{T_a\}$.
\end{proposition}

\begin{proof} 
  The first statement follows from
  Proposition~\ref{prop:sharpsections} and \eqref{eq:flatc_sections}.
  Since $[N_j,N_l]=0$ for all $1\leq j, l\leq r$, we have
  \begin{equation*}
    \begin{split}
      N_l(\ti{T}_a)\ &=\
      N_l\biggl(\exp\biggl(\sum_{j=1}^r\frac{\log(q_j)}{2\pi i}
      N_j\biggr) T_a^\flat\biggr)\ =\
      \exp\biggl(\sum_{j=1}^r\frac{\log q_j}{2\pi i} N_j \biggr)
      N_l(T_a^\flat)\\ &=\ \sum_b (L_{T_l})_{b a} \ti{T}_b, 
    \end{split}
  \end{equation*}
  and the second statement follows.
\end{proof}

\begin{remark}
  In the context of the Gromov-Witten potential, the previous result
  reduces to~\cite[Prop. 8.5.4]{bo:CK-mirror} whose proof involves the
  formalism of \jdef{gravitational correlators}.  The elementary proof
  given above shows that it is a direct consequence of the definition
  of the connection and, in particular, of the homogeneity of the the
  operators $L_{T_j}$.
\end{remark}

Because of Propositions~\ref{prop:sharpsections}
and~\ref{prop:david_854}, we know the first (graded) component of the
sections $T_a^\flat$ and $\ti{T}_a$. A lengthy but straightforward
computation yields the second component of both $T_a^\flat$ and
$\ti{T}_a$.

\begin{lemma}
  The $\nabla$-flat sections $T_a^\flat$ satisfy 
  \begin{equation*}
    T_a^\flat(q) \ \equiv\
    T_a - \sum_{\ti{c}=\ti{a}+2} \pdii{\phi}{a}{\delta(c)} T_c
    \mod \UU_{\ti{a}+2}.    
  \end{equation*}
\end{lemma}

\begin{lemma}\label{le:nablac_flat_sections}
  The $\nabla^c$-flat sections $\ti{T}_a$ satisfy the following
  formulas, for $k>3$. 

  For $\ti{a} \geq 2k-2$, $\ti{T}_a = T_a$.

  For $\ti{a}=2k-4$, $\ti{T}_a = T_a - \sum_{\ti{c}=\ti{a}+2} 2\pi i
  q_{\delta(c)} \pd{}{q_{\delta(c)}} \phi^a_h\, T_c + \phi^a_h\, T_m$.

  For $2<\ti{a}<2k-4$, $\ti{T}_a \equiv T_a - \sum_{\ti{c}=\ti{a}+2}
  \phi^{a \delta(c)}_h\, T_c \mod \UU_{\ti{a}+2}$.

  For $\ti{a}=2$, $\ti{T}_a \equiv T_a -
  \sum_{\ti{c}=\ti{a}+2} 2\pi i q_a \pd{}{q_a} \phi^{\delta(c)}_h\, T_c
  \mod \UU_{\ti{a}+2}$.

  For $\ti{a}=0$, $\ti{T}_0 \equiv T_0 \mod \UU_{\ti{a}+2}$.
\end{lemma}

We can now extend trivially the form $Q$, defined
by~\eqref{eq:Q_from_B_again}, to a form $\QQ$ on $\VV$. $\QQ$ is flat
because of~\eqref{eq:sym_condition}.  

To define a flat real structure $\VV_\R$ on $\VV$ we proceed as
follows. Let
\begin{equation*}
  \ti{\VV} \ := \ \Delta^r \times V \quad \text{ and } \quad \ti{\nabla} \
  := \ \nabla -\frac{1}{2\pi i} \sum_{j=1}^r N_j \frac{dq_j}{q_j}. 
\end{equation*}
Then $\ti{\nabla}$ is a flat connection on the bundle $\ti{\VV}$; for
$v \in V$ we define $\ti{\sigma}_v$ to be the $\ti{\nabla}$-flat
section of $\ti{\VV}$ such that $\ti{\sigma}_v(0) = v$. Then $\VV_\R$
is the local system generated by the sections $\exp(-\frac{1}{2\pi i}
\sum_{j=1}^r \log(q_j) N_j)\, \ti{\sigma}_v(q)$, for all $v\in V_\R$.

\begin{definition}
  \label{def:A-variation}
  Let $\phi = \phi_0 + \phi_\h$ be a quantum potential on the framed,
  real Frobenius module $(V,\CB,e,*)$.  Then $(
  \VV,\nabla,\FF,\VV_\R,\QQ)$ is the \jdef{$A$-model variation} of the
  potential.
\end{definition}

\begin{theorem}
  \label{th:a_model_is_pvhs}
  The $A$-model variation is a polarized variation of Hodge structure.
  Moreover, it is the variation associated to the potential $\phi$ by
  Theorem~\ref{th:correspondence_grl}.
\end{theorem}

\begin{proof}
  Let $\Phi$ be the ``period map'' of $( \VV,\nabla,\FF,\VV_\R,\QQ)$
  defined by parallel transport to the fiber $\VV_{\hat{q}}$,
  $\hat{q}\in(\Delta^*)^r$.  By Proposition~\ref{prop:N_product_0} the
  local monodromy logarithms $N_j$ are the left multiplication
  operators $L_{T_j}$ and, by Proposition~\ref{prop:david_854}, the
  limiting Hodge filtration becomes $F^p := \oplus_{a\geq p}
  V_{2(k-a)}$. Thus, Proposition~\ref{prop:framed_FM_NO} implies that
  $\{N_1,\dots,N_r;F\}$ is a nilpotent orbit.
  
  Let now $\exp (- \sum_j z_j\,N_j)\cdot \Phi(q) = \exp \Gamma(q)
  \cdot F$, where $\Gamma$ is a holomorphic, $\jgg_-$-valued map
  defined locally around $0\in \Delta^r$.  Since the map $\Phi$ is
  horizontal, the $\jgp_{-1}$-valued map $X_{-1} = \sum_j z_j\,N_j +
  \Gamma_{-1}$ satisfies the integrability condition
  \eqref{eq:integcond} and it follows from
  Theorem~\ref{th:improved_2.8} that $( \VV,\nabla,\FF,\VV_\R,\QQ)$ is
  a polarized variation of Hodge structure.
  
  In order to prove that this variation agrees with the one defined in
  Theorem~\ref{th:correspondence_grl} we appeal to the uniqueness
  statement in Theorem~\ref{th:improved_2.8}.  Hence, it suffices to
  show that $\Gamma_{-1}$ is related to the potential $\phi$ by
  \eqref{eq:gamma_from_potential}.  But, the matrix presentation of
  $\exp(-\Gamma(q))$ in the basis $\{T_a\}$ is the matrix expressing
  the $\nabla^c$-flat frame $\{\ti{T}_a\}$ in terms of the constant
  frame $\{T_a\}$. Thus, it follows from
  Lemma~\ref{le:nablac_flat_sections}, that
  \begin{equation*}
    \Gamma_{-1}(T_a)\ =\ \sum_{\ti{c}=\ti{a}+2}
    \pdii{\phi_\h}{a}{\delta(c)} T_c,
  \end{equation*}
  as desired.
\end{proof}


\def\cprime{$'$}
\providecommand{\bysame}{\leavevmode\hbox to3em{\hrulefill}\thinspace}


\end{document}